\documentclass[12pt]{article}
\usepackage{amsfonts,latexsym,amssymb}
 \newcommand{\beqn}{\begin{eqnarray}}
 \newcommand{\eeqn}{\end{eqnarray}}
 \newcommand{\be}{\begin{equation}}
 \newcommand{\ee}{\end{equation}}
 \newcommand{\ba}{\begin{array}}
 \newcommand{\ea}{\end{array}}
 \newcommand{\pa}{\partial}
 
 \newcommand{\ci}{\cite}
 \newcommand{\la}{\label}

\newcommand{\ga}{\gamma}

\newcommand{\na}{\nabla}

\newcommand{\bt}{\beta}
\newcommand{\al}{\alpha}

 \newcommand{\De}{\Delta}

\def\R{{\rm I\kern-.1567em R}}
\def\M{{\rm I\kern-.1567em M}}

\def\div {{\rm div}}

\def\spt{{\rm spt}}

 \newtheorem{theorem}{Theorem}[section]
 
 \newtheorem{definition}[theorem]{Definition}
 
 \newtheorem{lemma}[theorem]{Lemma}

 \newtheorem{pro}[theorem]{Proposition}

\begin{document}

\begin{center} {\Large\bf A note on necessary conditions for blow-up of energy
solutions to the Navier-Stokes equations}\\
  \vspace{1cm}
 {\large
  G. Seregin}

 \end{center}

 \vspace{1cm}
 \noindent
 {\bf Abstract } In the present note, we address the question about behavior
 of $L_3$-norm of the velocity field as time $t$ approaches blow-up time $T$.
 It is known that the upper limit of the above norm  must
 be equal to infinity. We show that, for blow-ups of type I, the lower limit of  $L_3$-norm  equals to infinity as well.

 \vspace {1cm}

\noindent {\bf 1991 Mathematical subject classification (Amer.
Math. Soc.)}: 35K, 76D.


\setcounter{equation}{0}
\section{Motivation }

Consider the Cauchy problem for the classical 3D-Navier-Stokes system

\be\la{m1}\left.
\begin{array}{c}\pa_t v+v\cdot\na\,v -\nu\,\Delta \,v =-\nabla
\,q  \\

            \div \,v = 0\end{array}
            \right\}\qquad\mbox{in}\,\,\, Q_+,\ee

\be\la{m2}v|_{t=0}=a\in C^\infty_{0,0}(\mathbb R^3).\ee Here, $v$
and $q$ stand for the velocity field and for the pressure field,
respectively, $Q_+=\mathbb R^3\times ]0,+\infty [$, and
$$C^\infty_{0,0}(\mathbb R^3)=\{a\in C^\infty_0(\mathbb R^3):\,\,\div\,a=0\,
\,\mbox{in}\,\,\mathbb R^3\,\}.$$
In what follows, we always assume that $\nu=1$.

It is well known due to J. Leray, see \ci{Le},  the Cauchy
problem (\ref{m1}), (\ref{m2}) has at least one solution called the
weak Leray-Hopf solution. To give its modern definition, let us
introduce standard energy spaces $H$ and $V$. $H$ is the closure of
the set $C^\infty_{0,0}(\mathbb R^3)$ in $L_2(\mathbb R^3)$ and $V$
is the closure of the same set with respect to the norm generated by
the Dirichlet integral.
\begin{definition}\la{md1} A velocity field $v\in L_\infty(0,+\infty;H)\cap L_2(0,+\infty;V)$
is called a weak Leray-Hopf solution to the Cauchy problem
(\ref{m1}),  (\ref{m2}) if the following conditions hold:
\be\la{m3}\int\limits_{Q_+}(v\cdot\pa_tw+v\otimes
v:\na\,w-\na\,v:\na\, w)\,dxdt=0\ee for any $w\in C_0^\infty(Q_+)$
with $\div\, w=0$ in $Q_+$;

the function
\be\la{m4}t\mapsto \int\limits_{\mathbb R^3}v(x,t)\cdot w(x)dx\ee
is continuous on $[0,+\infty [$ for all $w\in L_2(\mathbb R^3)$;
\be\la{m5}\|v(\cdot,t)-a(\cdot)\|_2\to 0\ee
as $t\to +0$;
\be\la{m6}\|v(\cdot,t)\|_2^2+2\int^t_0\int\limits_{\mathbb R^3}|\na\,v|^2dxdt'\leq
\|a\|_2^2\ee
for all $t\in [0,+\infty [$.
\end{definition}

The definition does not contain the pressure field at all. However,
using the linear theory, we can introduce so-called associated
pressure $q(\cdot,t)$, which,
for all $t>0$, 
is the Newtonian potential of
$v_{i,j}(\cdot,t)v_{j,i}(\cdot,t)$ and satisfies the pressure
equation \be\la{m7} -\De\,
q(\cdot,t)=v_{i,j}(\cdot,t)v_{j,i}(\cdot,t)=\div\,\div\,
v(\cdot,t)\otimes v(\cdot,t)\ee in $\mathbb R^3$. Since $v$ is known
to belong 
$L_\frac {10}3(Q_+)$, pressure
$q$
is in $L_\frac {5}3(Q_+)$. Moreover, the Navier-Stokes system is
satisfied in the sense of distributions and even
a.e. in $Q_+$.  We refer to the paper \ci{LS} for details.

Uniqueness of weak Leray-Hopf solutions is still unknown. However,
there is a simple but deep connection between smoothness and
uniqueness. It has been pointed out by J. Leray in his celebrated paper \ci{Le} and reads: any smooth
solution to (\ref{m1}),  (\ref{m2}) is unique in the class of weak
Leray-Hopf solutions. The problem of smoothness of weak Leray-Hopf
solutions is actually one of the seven Millennium problems.

In the paper, we
deal with certain necessary conditions for possible blow-ups of
solutions to the Cauchy problem (\ref{m1}),  (\ref{m2}). Suppose
that $T>0$ is the first moment of time when singularities occur.
Then, as it has been shown by J. Leray, given $3<s\leq +\infty$, there
exists a constant $c_s$ such that
\be\la{m8}\|v(\cdot,t)\|_{s,\mathbb R^3}\geq \frac {c_s}{(T-t)^\frac
{s-3}{2s}}\ee for $T/2\leq t< T$.

However, in the marginal case $s=3$, we have a weaker result
\be\la{m9}\limsup\limits_{t\to T-0}\|v(\cdot,t)\|_3=+\infty,\ee which has been established in
\ci{ESS4}. Apparently,  a natural question can be raised 
whether the  statement \be\la{m10}\lim\limits_{t\to
T-0}\|v(\cdot,t)\|_3=+\infty\ee  is true or not. In \ci{S5}, there
has been proved a weaker version of (\ref{m10}), namely,
\be\la{m11}\lim\limits_{t\to T-0}\frac
1{T-t}\int\limits^T_t\|v(\cdot,\tau)\|^3_3d\tau=+\infty.\ee

The aim of the present paper is to show validity of (\ref{m10})
provided  the blow-up of type I takes place, i.e.,
\be\la{m12}\|v(\cdot,t)\|_\infty\leq \frac {C_\infty}{\sqrt{T-t}}\ee
for any $T/2\leq t< T$ and for some positive constant $C_\infty$.
Our main result can be formulated as follows.
\begin{theorem}\la{mt2}
Let $T$ be a blow-up time and let, for some $3<s\leq +\infty$, there exist a positive
constant $C_s$ such that \be\la{m13}
\|v(\cdot,t)\|_{s}\leq \frac {C_s}{(T-t)^\frac
{s-3}{2s}}\ee for any $T/2\leq t< T$. Then (\ref{m10}) holds.
\end{theorem}

Let us outline our proof of
Theorem \ref{mt2}. 
 Firstly, we reduce the general
case to a particular one $s=+\infty$ showing that if (\ref{m13}) is
true for some $3<s<+\infty$, then it is true for $s=+\infty$ as
well. Secondly, assuming that (\ref{m10}) is violated, i.e., 
a sequence $t_k$ tending to 
$T$ exists such that 
\be\la{m14}\sup\limits_k \|v(\cdot,t_k)\|_3=M<+\infty,\ee we may use
a  blow-up machinery and construct a non-trivial ancient solution
defined in $\mathbb R^3\times ]-\infty,0[$ with the following
properties. It vanishes at time $t=0$ and  its $L_3$-norm is finite
say at time $t=-1$. In order to apply backward uniqueness results,
proved in \ci{ESS4}, we need
to check that the
above ancient solution has a certain behavior at infinity with respect to spatial variables. This can be done with the help of the
conception of so-called local energy solutions to the Cauchy
problem, see  \ci{LR1} and   also \ci{KS}.

Finally, it is interesting to figure out whether condition
(\ref{m14}) itself implies regularity. It is worthy to note that the
important consequence of (\ref{m14})
is that \be\la{m15}v(\cdot,T)\in L_3(\mathbb R^3).\ee

\setcounter{equation}{0}
\section{Some auxiliary things  }

In the paper, we are going to use the following notion. $B(x_0,R)$
stands for a spatial ball centered at a point $x_0$ and having
radius $R$, $B(R)=B(0,R)$, and $B=B(1)$. By $Q(z_0,R)$, where
$z_0=(x_0,t_0)$ is a space-time point, we denote a parabolic ball
$B(x_0,R)\times ]t_0-R^2,t_0[$, and $Q(R)=Q(0,R)$, $Q=Q(1)$. All
constants depending on non-essential parameters will be denoted
simply by $c$.
\begin{lemma}\la{al1} Suppose that (\ref{m13}) holds for some $3<s<+\infty$.
Then it is true for $s=+\infty$. \end{lemma} \textsc{Proof} From
(\ref{m7}) and (\ref{m13}) it follows that $q(\cdot,t)\in L_\frac
s2(\mathbb R^3)$ for $T/2<t<T$.

Fix $\varepsilon>0$ and $z_0=(x_0,t_0)$ with $t_0<T$ arbitrarily.  Applying (\ref{m13})
and H\"older inequality, we find
$$\frac 1{R^2}\int\limits_{Q(z_0,R)}(|v|^3+|q|^\frac 32)dz\leq$$$$\leq c(s)\frac 1{R^2}\Big(
\int\limits_{Q(z_0,R)}(|v|^s+|q|^\frac s2)dz\Big)^\frac 3sR^{5(1-\frac 3s)}\leq$$
$$\leq c(s)R^{5(1-\frac 3s)-2}\Big(
\int\limits^{t_0}_{t_0-R^2}\frac {C^s_s}{(T-t)^\frac {s-3}2}dt\Big)^\frac 3s\leq$$
$$\leq c(s)C^3_sR^{3-\frac {15}s}R^{2\frac 3s}\frac {1}{(T-t_0)^{3\frac {s-3}{2s}}}$$
$$\leq c(s)C^3_s\Big(\frac R{\sqrt{T-t_0}}\Big)^{3\frac {s-3}{2s}}.$$
We let $R=\sqrt{\ga(T-t_0)}$ and pick up $0<\ga<1$ so that
$c(s)C^3_s\ga^{3\frac {s-3}{2s}}\leq \varepsilon/2$. Now, we apply
the local regularity theory for suitable weak solutions to the
Navier-Stokes equations, developed in  \ci{CKN}, \ci{Li}, \ci{LS}, and
\ci{ESS4}. It reads that if $\varepsilon\leq \varepsilon_0$, where
$\varepsilon_0$ is a universal constant, then
$$|v(z_0)|\leq \frac {c}R= \frac {c}{\sqrt{\ga(T-t_0)}}$$ for all
$z_0=(x_0,t_0)$ with $x_0\in\mathbb R^3$ and $T/2<t_0<T$ and for some universal constant $c$.
Lemma \ref{al1} is proved.

So, we  need prove Theorem \ref{mt2} in a particular case $s=+\infty$ only.

\setcounter{equation}{0}
\section{Ancient Solution }

By assumptions of Theorem \ref{mt2}, there must be  singular points
at $t=T$. We take any of them say $x_0\in\mathbb R^3$. Then local
regularity theory gives the following inequality \be\la{s1} \frac
1{R^2}\int\limits_{Q((x_0,T),R)}(|v|^3+|q|^\frac 32)dz\geq
\varepsilon_0>0\ee for $0<R<R_0=\frac 13\min\{1,\sqrt{T}\}$ with
universal constant $\varepsilon_0$. Without loss of generality, we
may assume that $x_0=0$.

Proceeding in the same way as in \ci{SS2}, we can find that condition (\ref{m13}) implies the
following bound
$$\sup\limits_{0<R\leq R_0}\Big\{ \frac 1{R^2}\int\limits_{Q((0,T),R)}
(|v|^3+|q|^\frac 32)dz+\frac
1R\int\limits_{Q((0,T),R)}|\na\,v|^2dz+$$ \be\la{s2}+\frac
1R\sup\limits_{T-R^2<t<T}\int\limits_{B(R)}|v(x,t)|^2dx\Big\}=M_1<+\infty.\ee

Next, we may scale our functions $v$ and $q$ essentially in the same
way as it has been done in \ci{S5}, namely,
$$u^{(k)}(y,s)=R_kv(R_ky,T+R^2_ks),\qquad p^{(k)}(y,s)=R^2_kq(R_ky,T+R^2_ks)$$
for $y\in B(R_0/R_k)$ and for $s\in ]-(R_0/R_k)^2,0[$, where $R_k=\sqrt{T-t_k}$.

Now, let us see what happens if $k\to +\infty$. This is more or less
well-understood procedure and the reader can find  details 
in \ci{ESS4}, \ci{S5}--
\ci{SZ}. As a result, we
have two measurable functions $u$ and $p$ defined on $Q_-=\mathbb
R^3\times ]-\infty,0[$ with the following properties:
$$u^{(k)}\to u \qquad \mbox{in}\quad L_3(Q(a)),$$
$$\na\,u^{(k)}\rightharpoonup\na\, u \qquad \mbox{in}\quad L_2(Q(a)),$$
\be\la{s3}p^{(k)}\rightharpoonup p \qquad \mbox{in}\quad L_\frac 32(Q(a)),\ee
$$u^{(k)}\to u \qquad \mbox{in}\quad C([-a^2,0];L_\frac 98(B(a)))$$
for any $a>0$. The pair $u$ and $p$ satisfies the Navier-Stokes
equations in $Q_-$ in the sense distributions. We call  it  an
ancient solution to the Navier-Stokes equations. Moreover, since
inequalities (\ref{m12}) and (\ref{s2}) are invariant with respect
to the Navier-Stokes scaling, we can show that
$$\sup\limits_{0<a<+\infty}
\Big\{ \frac 1{a^2}\int\limits_{Q(a)}(|u|^3+|p|^\frac 32)de+\frac 1a\int\limits_{Q(a)}|\na\,u|^2de+$$
\be\la{s4}+\frac 1a\sup\limits_{-a^2<s<0}\int\limits_{B(a)}|u(y,s)|^2dy\Big\}\leq M_1<+\infty\ee
and
\be\la{s5}|u(y,s)|\leq \frac {C_\infty}{\sqrt{-s}}\ee
for all $e=(y,s)\in Q_-$.

The important consequence of (\ref{m15}) and
the last line in (\ref{s3}), is the following fact
\be\la{s6}u(\cdot,0)=0\ee in $\mathbb R^3$, see \ci{S5}
in a similar situation.

Now, our goal is to show that the above ancient solution is non-trivial. Unfortunately, we cannot get this by direct passing to the
limit in the formula
$$\frac 1{a^2}\int\limits_{Q(a)}(|u^{(k)}|^3+|p^{(k)}|^\frac 32)de=$$
\be\la{s7}=\frac 1{a^2R^2_k}\int\limits_{Q(aR_k)}(|v|^3+|q|^\frac 32)dz\geq \varepsilon_0>0\ee
for $aR_k<3/4$. The reason is simple: there is no hope to prove  strong convergence of the pressure. However, we still have local strong convergence of $u^{(k)}$ so that
\be\la{s8}\frac 1{a^2}\int\limits_{Q(a)}|u^{(k)}|^3de\to \frac 1{a^2}\int\limits_{Q(a)}|u|^3de\ee
for any $0<a\leq 3/4$.

To prove that our ancient solution is non-trivial, let us first note
that according to (\ref{s2}) \be\la{s9}\frac
1{a^2}\int\limits_{Q(a)}(|u^{(k)}|^3+|p^{(k)}|^\frac 32)de\leq
M_1\ee
 for sufficient large $k$   and  for all $a\in ]0,3/4]$.

The second observation is quite typical when treating the pressure.
In the ball $B(3/4)$, the pressure can be split into
two parts
$$p^{(k)}=p^{(k)}_1+p^{(k)}_2,$$
 where the first term is defined by the variational identity
$$\int\limits_{B(3/4)}p^{(k)}_1(y,s)\De\,\varphi(y)dy=-\int\limits_{B(3/4)}u^{(k)}(y,s)
\otimes u^{(k)}(y,s):\na^2\varphi(y)dy$$ being valid for any
$\varphi\in W^2_3(B(3/4))$ with $\varphi=0$ on $\pa B(3/4)$. It is
not difficult to show that the first counter-part of the pressure
satisfies the estimate \be\la{s10}\|p^{(k)}_1(\cdot,s)\|_{\frac
32,B(3/4)}\leq c \|u^{(k)}(\cdot,s)\|^2_{3,B(3/4)}\ee for all
$-\infty <s<0$ while the second one
 is a harmonic function in $B(3/4)$
for the same $s$. Since $p^{(k)}_2(\cdot,s)$ is harmonic, we have
$$\sup\limits_{y\in B(1/2)}|p^{(k)}_2(y,s)|^\frac 32\leq c \int
\limits_{B(3/4)}|p^{(k)}_2(y,s)|^\frac 32dy\leq$$
\be\la{s11}\leq c\int\limits_{B(3/4)}|p^{(k)}(y,s)|^\frac 32dy+c\int
\limits_{B(3/4)}|u^{(k)}(y,s)|^ 3dy. \ee
Then, for any $0<a<1/2$,
$$\varepsilon_0\leq \frac 1{a^2}\int\limits_{Q(a)}(|u^{(k)}|^3+|p^{(k)}|^\frac 32)de
\leq$$$$\leq c\frac 1{a^2}\int\limits_{Q(a)}(|u^{(k)}|^3+|p^{(k)}_1|^\frac 32
+|p^{(k)}_2|^\frac 32)de$$$$\leq c\frac 1{a^2}\int\limits_{Q(a)}(|u^{(k)}|^3
+|p^{(k)}_1|^\frac 32)de+$$$$
+ca^3\frac 1{a^2}\int\limits_{-a^2}^0\sup\limits_{y\in B(1/2)}|p^{(k)}_2(y,s)|^\frac 32ds.$$
Combining  (\ref{s9})--(\ref{s11}),  we find
$$\varepsilon_0\leq 
c\frac 1{a^2}\int\limits_{Q(3/4)} |u^{(k)}|^ 3de +ca\int\limits_{-a^2}^0ds\int\limits_{B(3/4)}\Big(|p^{(k)}(y,s)|^\frac 32+
|u^{(k)}(y,s)|^ 3\Big)dy\leq$$
$$
\leq c\frac 1{a^2}\int\limits_{Q(3/4)} |u^{(k)}|^ 3de +ca\int\limits_{Q(3/4)}\Big(|p^{(k)}|^\frac 32+
|u^{(k)}|^ 3\Big)de \leq $$
$$\leq c\frac 1{a^2}\int\limits_{Q(3/4)} |u^{(k)}|^ 3de+cM_1a$$
for the same $a$.
 Passing to the limit and choosing sufficiently small $a$,
we show
that
\be\la{s12}0<c\varepsilon_0a^2\leq\int\limits_{Q(3/4)} |u|^ 3de\ee
for some positive $0<a<1/2$. So, our ancient solution $u$ is non-trivial.

If would show that for some positive $R_*$
$$|u|+|\na\,u|\in L_\infty((\mathbb R^3\setminus B(R_*))\times]-(5/6)^2,0[),$$
we could use arguments from \ci{ESS4} and conclude that, by
(\ref{s6}), $\na\wedge u\equiv 0$ in $\mathbb R^3\times
]-(3/4)^2,0[$ which, together with the incompressibility condition,
means that $u(\cdot,t)$ is harmonic in $\mathbb R^3$. And it is bounded there.
So, $u$ must be a function of $t$ only. But estimate (\ref{s4}) says
that such a function must be zero in $]-(3/4)^2,0[$. The latter
contradicts (\ref{s12}).

\setcounter{equation}{0}

\section{Spatial decay for ancient solutions }

We know that
$$\|u^{(k)}(\cdot,-1)\|_3\leq M$$
and thus by (\ref{s3})
\be\la{d1}\|u(\cdot,-1)\|_3\leq M.\ee

Now, let us consider 
the following Cauchy problem

\be\la{d2}\left.
\begin{array}{c}\pa_t w+w\cdot\na\,w -\Delta \,w =-\nabla
\,r  \\

            \div \,w = 0\end{array}
            \right\}\qquad\mbox{in}\,\,\,\widetilde{Q}=\mathbb R^3\times ]-1,1[,\ee

\be\la{d3}w(\cdot,-1)=u(\cdot,-1).\ee
We would like to construct a solution to problem (\ref{d2}), (\ref{d3})
satisfying the local energy inequality.
To this end, 
let us recall notation and some facts from \ci{KS}.
$$ L_{m,unif}=\{u\in L_{m,loc} \,\,:\,\,\|u\|_{L_{m,unif}}=
\sup\limits_{x_0\in\mathbb
R^3}\Big(\int\limits_{B(x_0,1)}|u(x)|^mdx\Big)^{1/m}< +\infty\},$$
$$E_m=\{u\in L_{m,unif}\,\,:\,\,\int\limits_{B(x_0,1)}|u(x)|^mdx\to
0\quad\mbox{as}\quad |x_0|\to +\infty\},$$
$$\stackrel{\circ}{E}_m=\{u\in E_m\,\,:\,\,\div\, u=0\quad
\mbox{in}\quad \mathbb R^3\}.$$
Apparently,
\be\la{d4}u(\cdot,-1)\in\stackrel{\circ}{E}_2.\ee
\begin{definition}\la{dd1}
A pair of functions $w$ and $r$ defined in the space-time
cylinder $\widetilde{Q}$ is called a local energy weak
 Leray-Hopf solution or simply local energy solution
 to the Cauchy problem (\ref{d2}), (\ref{d3}) if
 the following conditions are satisfied:
 $$
w\in L_\infty(-1,1;L_{2,unif}),\qquad\sup\limits_{x_0\in \mathbb R^3}\int\limits^1_{-1}
\int\limits_{B(x_0,1)}|\na w|^2dz<+\infty, $$
\be\la{d5} r\in L_\frac
32(-1,1;L_{\frac 32,loc}(\mathbb R^3));\ee
\be\la{d6}w\,\, and\,\,r\,\, meet\,\, (\ref{d2})\,\,in
\,\,the\,\,sense\,\, of\,\,distributions;\ee \be\la{d7} the\,
function\,\,t\mapsto\int\limits_{\mathbb R^3}w(x,t)\cdot \widetilde{w}(x)\,dx
\,is\,\, continuous\,\,on\,\,[-1,1] \ee for any compactly supported
function $\widetilde{w}\in L_2(\mathbb R^3)$;

for any compact K, \be\la{d8}\|w(\cdot,t)-u(\cdot,-1)\|_{L_{2}(K)}\to 0
\quad as\quad t\to -1+0;\ee

$$\nonumber \int\limits_{\mathbb R^3}\varphi|w(x,t)|^2\,dx+
2\int\limits_{-1}^{t}\int\limits_{\mathbb R^3}\varphi|\na
w|^2\,dxdt\leq \int\limits_{-1}^{t}\int\limits_{\mathbb
R^3}\Big(|w|^2(\pa_t\varphi+\Delta\varphi) +$$\begin{eqnarray}\la{d9}
 +w\cdot\na \varphi(|w|^2+2r)\Big)\,dxdt  &&
\end{eqnarray}
for a.a. $t\in ]-1,1[$ and for all nonnegative  functions
$\varphi\in C^\infty_0(\mathbb R^3\times ]-1,2[)$;

 for any
$x_0\in\mathbb R^3$, there exists a function $c_{x_0}\in L_\frac 32
(-1,1)$ such that \be\la{d10}r_{x_0}(x,t)\equiv
r(x,t)-c_{x_0}(t)=r_{x_0}^1(x,t)+r_{x_0}^2(x,t),\ee for $(x,t)\in
B(x_0,3/2)\times ]-1,1[$, where
$$r_{x_0}^1(x,t)=-\frac 13 |w(x,t)|^2
+\frac 1{4\pi}\int\limits_{B(x_0,2)}K(x-y): w(y,t)\otimes
w(y,t)\,dy,$$$$r_{x_0}^2(x,t)=\frac 1{4\pi}\int\limits_{\mathbb
R^3\setminus B(x_0,2)}(K(x-y)-K(x_0-y)): w(y,t)\otimes w(y,t)\,dy$$
and $K(x)=\na^2 (1/|x|)$.
\end{definition}
We have, see \ci{LR1} and also  \ci{KS}.
\begin{pro}\la{dp2}Under assumption (\ref{d4}), there exists at least
one local energy solution to problem (\ref{d2}), (\ref{d3}). \end{pro}
To describe spatial decay of local energy solution, we need additional notation
$$\alpha_w(t)=\|w(\cdot,t)\|^2_{L_{2,unif}},\qquad
\beta_w(t)=\sup\limits_{x_0\in \mathbb R^3}\int\limits^t_{-1}
\int\limits_{B(x_0,1)}|\na w|^2dxdt',$$
$$\ga_w(t)=\sup\limits_{x_0\in \mathbb R^3}\int\limits^t_{-1}
\int\limits_{B(x_0,1)}|w|^3dxdt',\qquad \delta_r(t)=
\sup\limits_{x_0\in \mathbb R^3}\int\limits^t_{-1}
\int\limits_{B(x_0,3/2)}|r_{x_0}|^\frac 32dxdt'.$$ One of the most
important properties of local energy solutions is a kind of uniform local boundedness of the energy, i.e.,
\be\la{d11}\sup\limits_{-1\leq t\leq
1}\alpha_w(t)+\beta_w(1)+\ga^\frac 23_w(1)\leq A<+\infty.\ee

Next, fix a smooth cut-off function $\chi$ so that $\chi(x)=0$
if $x\in B$, $\chi(x)=1$ if $x\notin B(2)$, and then let $\chi_R(x)=\chi(x/R)$.
Hence, one can define
$$\alpha^R_w(t)=\|\chi_Rw(\cdot,t)\|^2_{L_{2,unif}},
\qquad\beta^R_w(t)=\sup\limits_{x_0\in \mathbb R^3}\int\limits^t_{-1}
\int\limits_{B(x_0,1)}|\chi_R\na w|^2dxdt',$$
$$\ga^R_w(t)=\sup\limits_{x_0\in \mathbb R^3}\int\limits^t_{-1}
\int\limits_{B(x_0,1)}|\chi_R w|^3dxdt',\quad
\delta^R_r(t)=\sup\limits_{x_0\in \mathbb R^3}\int\limits^t_{-1}
\int\limits_{B(x_0,3/2)}|\chi_R r_{x_0}|^\frac 32dxdt'.$$
As it was shown in \ci{KS}, the following decay estimate is true.
\begin{lemma}\la{dl3} Assume that the pair $w$ and $r$ is a local
energy solution to (\ref{d2}), (\ref{d3}). Then
$$\sup\limits_{-1\leq t\leq }\alpha^R_w(t)+\beta^R_w(1)
+(\ga^R_w)^\frac 23(1)+(\delta^R_r)^\frac 43(1)\leq $$\be\la{d12}
\leq C(A)\Big[\|\chi_Ru(\cdot,-1)\|^2_{L_{2,unif}}+1/R^\frac 23\Big].\ee
\end{lemma}

Since any local energy solution to the Cauchy problem (\ref{d2}), (\ref{d3}) is also
a suitable weak solution to the Navier-Stokes equations, one can apply the local
regularity theory to them and deduce from Lemma \ref{dl3} that there exists a positive number $R_*$ such that
\be\la{d13}|w(z)|+|\na\,w(z)|\leq A_1\ee
for all $z=(x,t)\in \Big(\mathbb R^3\setminus B(R_*)\Big)\times [-(5/6)^2, 1].$

If we would show that
\be\la{d14}u\equiv w\ee
on $\mathbb R^3\times [0,1[$, this would make it possible to apply backward
uniqueness results (actually, to vorticity equations) and conclude that $u=0$
on $\mathbb R^3\times [-(3/4)^2, 0]$ which contradicts (\ref{s12}). So, the rest
of the paper is devoted to a proof of (\ref{d14}).

Our first observation in this direction is that $u$ is $C^\infty$-
function in $Q_-$. This follows from (\ref{s5}). Detail discussion  on differentiability properties of bounded ancient solutions can be found in
\ci{KNSS} and \ci{SS2}. In addition, the pressure $p(\cdot,t)$ is a BMO-solution to the pressure equations
$$-\De p(\cdot,t)=\div\div\, u(\cdot,t)\otimes u(\cdot,t)$$
in $R^3$.

Using a suitable cut-off function in time and differentiability
properties of $w$ and $u$, we can get the following three relations:
$$\int\limits_{\mathbb R^3}\varphi(x)w(x,\tau)\cdot u(\cdot,t)dx
\Big|^{\tau=t}_{\tau=-1}=$$$$=\int\limits_{-1}^t\int\limits_{\mathbb
R^3}(w\otimes w-\na\,w):(\varphi\na\,u+u\otimes
\na\,\varphi)dxd\tau+$$$$+\int\limits_{-1}^t\int\limits_{\mathbb
R^3}ru\cdot \na\,\varphi
dxd\tau+\int\limits_{-1}^t\int\limits_{\mathbb R^3} \phi
w\cdot\pa_tudxd\tau;$$
$$\int\limits_{\mathbb R^3}\varphi(x)|w(x,\tau)|^2dx
\Big|^{\tau=t}_{\tau=-1}+2\int\limits_{-1}^t\int\limits_{\mathbb
R^3}\varphi
|\na\,w|^2dxd\tau\leq$$$$\leq\int\limits_{-1}^t\int\limits_{\mathbb
R^3}\Big(|w|^2\De\,\varphi+w\cdot
\na\,\varphi(|w|^2+2r)\Big)dxd\tau;$$
$$\int\limits_{\mathbb R^3}\varphi(x)|u(x,\tau)|^2dx
\Big|^{\tau=t}_{\tau=-1}+2\int\limits_{-1}^t\int\limits_{\mathbb
R^3}\varphi
|\na\,u|^2dxd\tau=$$$$=\int\limits_{-1}^t\int\limits_{\mathbb
R^3}\Big(|u|^2\De\,\varphi+u\cdot
\na\,\varphi(|u|^2+2p)\Big)dxd\tau$$ for any $0\leq\varphi\in
C^\infty_0(\mathbb R^3)$. Letting $\overline{u}=w-u$ and
$\overline{p}=r-p$, we can find from them the main inequality
$$\int\limits_{\mathbb R^3}\varphi(x)|\overline{u}(x,t)|^2dx
+2\int\limits_{-1}^t\int\limits_{\mathbb R^3}\varphi
|\na\,\overline{u}|^2dxd\tau\leq$$\be\la{d15}\leq\int\limits_{-1}^t
\int\limits_{\mathbb
R^3}\Big(|\overline{u}|^2\De\,\varphi+\overline{u}\cdot
\na\,\varphi(|\overline{u}|^2+2\overline{p})+u\cdot
\na\,\varphi|\overline{u}|^2+\ee$$
-2\varphi\na\,u:\overline{u}\otimes\overline{u}\Big)dxd\tau$$ for
a.a. $t$ in $]-1,0[$.

Next, for $u$ and $\overline{u}$, we may introduce the analogous quantities
$$\alpha_u(t)=\|u(\cdot,t)\|^2_{L_{2,unif}},
\qquad\alpha_{\overline{u}}(t)=\|\overline{u}(\cdot,t)\|^2_{L_{2,unif}}$$

$$\beta_u(t)=\sup\limits_{x_0\in \mathbb R^3}\int\limits^t_{-1}
\int\limits_{B(x_0,1)}|\na\,
u|^2dxdt',\quad\beta_{\overline{u}}(t)=\sup\limits_{x_0\in \mathbb
R^3}\int\limits^t_{-1} \int\limits_{B(x_0,1)}|\na\,
\overline{u}|^2dxdt',$$
$$\ga_u(t)=\sup\limits_{x_0\in \mathbb R^3}\int\limits^t_{-1}
\int\limits_{B(x_0,1)}|
u|^3dxdt',\quad\ga_{\overline{u}}(t)=\sup\limits_{x_0\in \mathbb
R^3}\int\limits^t_{-1}
\int\limits_{B(x_0,1)}|\overline{u}|^3dxdt',$$$$
\delta_p(t)=\sup\limits_{x_0\in \mathbb R^3}\int\limits^t_{-1}
\int\limits_{B(x_0,3/2)}|p_{x_0}|^\frac
32dxdt',\quad\delta_{\overline{p}}(t)=\sup\limits_{x_0\in \mathbb
R^3}\int\limits^t_{-1}
\int\limits_{B(x_0,3/2)}|\overline{p}_{x_0}|^\frac 32dxdt'$$ where
$$p_{x_0}(x,t)\equiv
p(x,t)-p^{0}_{x_0}(t)=p_{x_0}^1(x,t)+p_{x_0}^2(x,t),$$
$$p_{x_0}^1(x,t)=-\frac 13 |u(x,t)|^2
+\frac 1{4\pi}\int\limits_{B(x_0,2)}K(x-y): u(y,t)\otimes
u(y,t)\,dy,$$$$p_{x_0}^2(x,t)=\frac 1{4\pi}\int\limits_{\mathbb
R^3\setminus B(x_0,2)}(K(x-y)-K(x_0-y)): u(y,t)\otimes u(y,t)\,dy,$$
$$\overline{p}_{x_0}(x,t)\equiv
\overline{p}(x,t)-\overline{p}^{0}_{x_0}(t)=
\overline{p}_{x_0}^1(x,t)+\overline{p}_{x_0}^2(x,t),$$
$$\overline{p}_{x_0}^1(x,t)=r_{x_0}^1(x,t)-p_{x_0}^1(x,t),$$
$$\overline{p}_{x_0}^2(x,t)=r_{x_0}^2(x,t)-p_{x_0}^2(x,t).$$

By (\ref{s5}) and by our definitions,
\be\la{d16}\alpha_u(t)\leq\frac c{(-t)},\quad \beta_u(t)\leq\frac
c{(-t)^2},\quad \ga_u(t)+\delta_p(t)\leq \frac c{(-t)^\frac 32}\ee
for all $-1\leq t<0$. Indeed, the first bound follows directly from (\ref{s5}).
To get the second one, we need to use (\ref{s5}), BMO-estimate of the pressure
via velocity field, and then local regularity theory in the same way as in the proof
of Lemma \ref{al1}. It is useful to note that the above arguments imply the estimate $|\na\,u(x,t)|\leq c/(-t)$ for all $(x,t)\in Q_-$. As to the third bound, the second term is estimated with the help
of the singular integral theory, (\ref{s5}), and definitions of $p_{x_0}^1$ and $p_{x_0}^2$.

We fix $x_0\in \mathbb R^3$ and a smooth non-negative function
$\varphi$ such that $\varphi\equiv 1$ in $B$ and
$\spt\,\varphi\subset B(3/2)$ and let
$\varphi_{x_0}(x)=\varphi(x-x_0)$. Considering  (\ref{d15}) with such a cut-off function
$\varphi_{x_0}$, taking into account (\ref{d11}) and
(\ref{d16}), and arguing for example as in \ci{KS}, we can find the
inequality
$$\alpha_{\overline{u}}(t_0)+\beta_{\overline{u}}(t_0)\leq
c\Big[\int\limits^{t_0}_{-1}\alpha_{\overline{u}}(t)dt+\ga_{\overline{u}}(t_0)+$$
\be\la{d17}+\sup\limits_{x_0\in \mathbb
R^3}\int\limits^{t_0}_{-1}\int\limits_{B(x_0,3/2)}|\overline{p}_{x_0}|
|\overline{u}|dxdt+\ee
$$+\frac
c{\sqrt{-t_0}}\int\limits^{t_0}_{-1}\alpha_{\overline{u}}(t)dt+
\frac c{(-t_0)}\int\limits^{t_0}_{-1}\alpha_{\overline{u}}(t)dt\Big]$$
for all $-1\leq t_0<0$.

Next, we can re-write the well-known (in the theory of the Navier-Stokes equations)  multiplicative inequality in terms of quantities introduced above
$$\ga_{\overline{u}}(t_0)\leq
c\Big(\int\limits^{t_0}_{-1}\alpha^3_{\overline{u}}(t)dt\Big)^\frac
14\Big(\beta_{\overline{u}}(t_0)+\int\limits^{t_0}_{-1}\alpha_{\overline{u}}(t)dt\Big)^\frac
34.$$ To simplify the latter, we first make use of (\ref{d11}) and
(\ref{d16}) in the following way
$$\alpha_{\overline{u}}(t_0)\leq
c(\alpha_{w}(t_0)+\alpha_{u}(t_0))\leq c\Big(A+\frac
c{(-t_0)}\Big)\leq \frac {C(A)}{(-t_0)}$$ for all $-1\leq t_0<0$. And
thus \be\la{d18}\ga_{\overline{u}}(t_0)\leq \frac
{C(A)}{\sqrt{-t_0}}
\Big(\int\limits^{t_0}_{-1}\alpha_{\overline{u}}(t)dt\Big)^\frac
14\Big(\beta_{\overline{u}}(t_0)+\int\limits^{t_0}_{-1}
\alpha_{\overline{u}}(t)dt\Big)^\frac
34\ee

It remains to estimate the third term on the right hand side  of (\ref{d17})
\be\la{d19}I=\sup\limits_{x_0\in \mathbb
R^3}\int\limits^{t_0}_{-1}\int\limits_{B(x_0,3/2)}
|\overline{p}_{x_0}||\overline{u}|dxdt\leq I_1+I_2,\ee
where
$$I_1=\sup\limits_{x_0\in \mathbb
R^3}\int\limits^{t_0}_{-1}\int\limits_{B(x_0,3/2)}
|\overline{p}^1_{x_0}||\overline{u}|dxdt\leq I'_1+I''_1,$$
$$I_2=\sup\limits_{x_0\in \mathbb
R^3}\int\limits^{t_0}_{-1}\int\limits_{B(x_0,3/2)}
|\overline{p}^2_{x_0}||\overline{u}|dxdt,$$
and
$$I'_1=c\sup\limits_{x_0\in \mathbb
R^3}\int\limits^{t_0}_{-1}\int\limits_{B(x_0,3/2)}\Big||w|^2-|u|^2\Big|
|\overline{u}|dxdt, $$
$$I''_1=c\sup\limits_{x_0\in \mathbb
R^3}\int\limits^{t_0}_{-1}\int\limits_{B(x_0,3/2)}\Big|\int\limits_{B(x_0,2)}
K(x-y):\Big(w(y,t)\otimes w(y,t)-$$$$-u(y,t)\otimes u(y,t)\Big)dy\Big||\overline{u}(x,t)|dxdt.$$
$I'_1$ is evaluated easily, namely,
$$I'_1\leq c\sup\limits_{x_0\in \mathbb
R^3}\int\limits^{t_0}_{-1}\int\limits_{B(x_0,3/2)}|\overline{u}|^2(|\overline{u}|
+2|u|)dxdt\leq $$
\be\la{d20}\leq c \ga_{\overline{u}}(t_0)+\frac c{\sqrt{-t_0}}\int\limits^{t_0}_{-1}\al_{\overline{u}}(t)dt.\ee
To estimate $I''_1$, we exploit the same idea and $L_\frac 32$- and $L_2$-estimates
for singular integrals
$$I''_1\leq c\sup\limits_{x_0\in \mathbb
R^3}\int\limits^{t_0}_{-1}\int\limits_{B(x_0,3/2)}|\overline{u}(x,t)|
\Big\{\Big|\int\limits_{B(x_0,2)}K(x-y):\overline{u}(y,t)\otimes \overline{u}(y,t)dy\Big|+$$$$+\Big|\int\limits_{B(x_0,2)}K(x-y):
\Big(\overline{u}(y,t)\otimes u(y,t)+u(y,t)\otimes \overline{u}(y,t)\Big)dy\Big|\Big\}dxdt\leq$$
$$\leq c \ga_{\overline{u}}(t_0)+\frac c{\sqrt{-t_0}}\int\limits^{t_0}_{-1}\al_{\overline{u}}(t)dt.$$
So, by (\ref{d20}), we have
\be\la{d21}I_1\leq c \ga_{\overline{u}}(t_0)+\frac c{\sqrt{-t_0}}\int\limits^{t_0}_{-1}\al_{\overline{u}}(t)dt.\ee

In order to find  upper bound for $I_2$, we simply repeat arguments of Lemma 2.1
in \ci{KS} with $R=1$ there. This gives us  the following estimate
$$|p^2_{x_0}(x,t)|\leq c\int\limits_{\mathbb R^3\setminus B(x_0,2)}\Big|K(x-y)-K(x_0-y)\Big|\Big|w(y,t)\otimes w(y,t)-$$$$-u(y,t)\otimes u(y,t)\Big|dx\leq$$
$$\leq c\|w(\cdot,t)\otimes w(\cdot,t)-u(\cdot,t)\otimes u(\cdot,t)\|_{L_{1,unif}}$$
being valid for any $x\in B(x_0,32/)$ and thus
$$I_2\leq c\sup\limits_{x_0\in \mathbb
R^3}\int\limits^{t_0}_{-1}\|w(\cdot,t)\otimes w(\cdot,t)-u(\cdot,t)\otimes u(\cdot,t)\|_{L_{1,unif}}\int\limits_{B(x_0,3/2)}|\overline{u}(x,t)|dx.$$
Furthermore, by (\ref{d11}),
$$\|w(\cdot,t)\otimes w(\cdot,t)-u(\cdot,t)\otimes u(\cdot,t)\|_{L_{1,unif}}=$$
$$=\sup\limits_{x_0\in \mathbb
R^3}\int\limits_{B(x_0,1)}|\overline{u}(y,t)\otimes w(y,t)+u(y,t)\otimes \overline{u}(y,t)|dx\leq$$
$$\leq \al_{\overline{u}}^\frac 12(t)\al_{w}^\frac 12(t)+\al_{\overline{u}}^\frac 12(t)\al_{u}^\frac 12(t)\leq \frac {C(A)}{\sqrt{-t}} \al_{\overline{u}}^\frac 12(t).$$
So,
$$I_2\leq \frac {C(A)}{\sqrt{-t_0}} \int\limits^{t_0}_{-1} \al_{\overline{u}}(t)dt$$
and, by (\ref{d18}) and (\ref{d21}), we have
$$I\leq c\ga_{\overline{u}}(t_0)+\frac {C(A)}{\sqrt{-t_0}}\int\limits^{t_0}_{-1}\al_{\overline{u}}(t)dt\leq$$
\be\la{d22}\leq \frac {C(A)}{\sqrt{-t_0}}\Big[\Big(\int\limits^{t_0}_{-1}\al_{\overline{u}}(t)dt\Big)^\frac 14\Big(\bt_{\overline{u}}(t_0)+\int\limits^{t_0}_{-1}
\al_{\overline{u}}(t)dt\Big)^\frac 34+\int\limits^{t_0}_{-1}\al_{\overline{u}}(t)dt\Big].\ee
Combining (\ref{d17}) and (\ref{d22}) and  applying Young inequality, we arrive at
the final estimate
$$\al_{\overline{u}}(t_0)\leq C(A,\delta)\int\limits^{t_0}_{-1}\al_{\overline{u}}(t)dt$$
 which is valid for all $-1\leq t_0\leq \delta<0$. The latter says that $\al_{\overline{u}}(t)=0$
in $[-1,0[$ and, hence, $u(\cdot,t)=w(\cdot,t)$ for the same $t$. This completes the proof of Theorem \ref{mt2}.

\noindent\textbf{Acknowledgement} This work was partially supported
by  the RFFI grant 08-01-00372-a.

\noindent

G. Seregin\\
Center for Nonlinear PDE's,\\
Mathematical Institute, University of Oxford,UK\\

\end{document}